\documentclass[11pt]{amsart}

\usepackage{amsmath}
\usepackage{amssymb}
\usepackage{sty1}
\usepackage{accents}

\makeatletter
\@addtoreset{equation}{section}
\makeatother

\def\am#1{\mathrm{AM}(#1)}
\def\bloneg{\mathrm{L}^1(G)}
\def\blonegn{\mathrm{L}^1(G/N)}
\def\blonegg{\mathrm{L}^1(G\cross G)}
\def\blone#1{\mathrm{L}^1(#1)}
\def\conjg{\mathrm{Conj}(G)}

\def\zbloneg{\mathrm{ZL}^1(G)}
\def\zblonegg{\mathrm{ZL}^1(G\cross G)}
\def\zblonegn{\mathrm{ZL}^1(G/N)}
\def\zblone#1{\mathrm{ZL}^1(#1)}

\def\blonek{\mathrm{L}^1(K)}
\def\conjsutwo{\mathrm{Conj}(\mathrm{SU}(2))}
\def\sutwo{\mathrm{SU}(2)}

\def\measg{\mathrm{M}(G)} 
\def\measgg{\mathrm{M}(G\cross G)}
\def\zmeasg{\mathrm{ZM}(G)} 
\def\zmeasgg{\mathrm{ZM}(G\cross G)} 

\def\z#1{\mathrm{Z}(#1)}

\def\trig{\mathrm{Trig}(G)}
\def\trigg{\mathrm{Trig}(G\cross G)}
\def\ztrig{\mathrm{ZTrig}(G)}
\def\zcontg{\mathrm{Z}\mathcal{C}(G)}
\def\zcontcg{\mathrm{Z}\mathcal{C}_c(G)}

\def\autg{\mathrm{Aut}(G)}
\def\inng{\mathrm{Inn}(G)}
\def\zbbloneg{\mathrm{Z}_B\mathrm{L}^1(G)}
\def\auth{\mathrm{Aut}(H)}
\def\innh{\mathrm{Inn}(H)}
\def\zbbloneh{\mathrm{Z}_B\mathrm{L}^1(H)}
\def\zbblonef{\mathrm{Z}_B\mathrm{L}^1(F)}

\def\chat#1{\accentset{\vspace{-0.6mm}\smallfrown}{#1}}

\begin{document}

\newtheorem{mainconj}{Conjecture}[section]
\newtheorem{hypertaubprop}[mainconj]{Theorem}

\newtheorem{tensorprod}{Proposition}[section]
\newtheorem{quotient}[tensorprod]{Proposition}
\newtheorem{quotient1}[tensorprod]{Corollary}
\newtheorem{weakstarclosed}[tensorprod]{Proposition}
\newtheorem{nicebai}[tensorprod]{Lemma}
\newtheorem{diagonal}[tensorprod]{Theorem}
\newtheorem{conngroup}[tensorprod]{Theorem}
\newtheorem{fingroup}[tensorprod]{Theorem}
\newtheorem{fingroup1}[tensorprod]{Corollary}
\newtheorem{prodfingroup}[tensorprod]{Theorem}
\newtheorem{finabexample}[tensorprod]{Theorem}
\newtheorem{sutwona}[tensorprod]{Theorem}
\newtheorem{teespfa}[tensorprod]{Theorem}

\newtheorem{stegmeir}{Lemma}[section]
\newtheorem{finitederived}[stegmeir]{Theorem}
\newtheorem{hypertauberian1}[stegmeir]{Proposition}
\newtheorem{hypertauberian2}[stegmeir]{Theorem}

\title[Centres of group algebras]
{Amenability properties of the centres of group algebras}

\author{Ahmadreza Azimifard, Ebrahim Samei and Nico Spronk}

\begin{abstract}
Let $G$ be a locally compact group, and $\zbloneg$ be the centre
of its group algebra.  We show that when $G$ is compact $\zbloneg$
is not amenable when $G$ is either nonabelian and connected, or is a 
product of infinitely many finite nonabelian groups.  We also, study,
for some non-compact groups $G$, some conditions which imply amenability
and hyper-Tauberian property, for $\zbloneg$.
\end{abstract}

\maketitle

\footnote{{\it Date}: \today.


Research of the first named author supported by the Workshop on Operator
Algebras (2007), at the Fields Institute.
Research of the second named author supported by an NSERC Post Doctoral 
Fellowship.  
Research of the the third named author supported by NSERC Grant 312515-05.}

Let $G$ be a locally compact group and $\bloneg$ denote the group
algebra, i.e.\ the subalgebra of the measure algebra $\measg$
consisting of measures which are absolutely continuous with respect to the
left Haar measure.  We let
\[
\zbloneg=\{f\in\bloneg:f\con g=g\con f\text{ for all }g\in\bloneg\}
\]
be the centre of $\bloneg$.  Our goal is to study amenability
and weak amenability for $\zbloneg$.

We show that when $G$ is compact, $\zbloneg$ is generally not amenable.
In fact, it fails to be amenable whenever $G$ is either non-abelian
and connected (Section \ref{ssec:conngroup}), or when $G$ is a product
of infinitely many non-abelian finite groups (Section \ref{ssec:prodgroup}).
These results make substantial use of some discoveries from the intensive study
of central idempotent measures on compact groups of D. Rider \cite{rider}.
It is mentioned in \cite{stegmeir} that it is known to B. E. Johnson that
$\zbloneg$ fails to be amenable for some compact group $G$.
However no further information is provided.
Our results, and techniques therein, lead us towards 
the following.

\begin{mainconj}\label{conj:mainconj}
If $G$ is compact, then $\zbloneg$ is amenable if and only
if $G$ admits an open abelian subgroup.
\end{mainconj}

\noindent
We note that by \cite{moore}, $G$ admits an open abelian subgroup if and
only if the set of degrees of its irreducible representations is bounded.
We address the above conjecture with an illustrative example (Section 
\ref{ssec:amenex}). 
As a complement to many of the methods used in the prior sections, we illustrate
two examples using hypergroup techniques (Section \ref{ssec:hypergroup}).

We close the article with a study of some non-compact groups.
We use results of R.D. Mosak and J. Liukkonen
\cite{mosak0,liukkonenm,mosak} extensively.
When the commutator of $G$ with the open subgroup, which
supports all elements of $\zbloneg$, is finite, then $\zbloneg$
is amenable (Section \ref{ssec:fincommutator}).  When $G$ has 
relatively compact conjugacy classes, then $\zbloneg$ is hyper-Tauberian
(Section \ref{ssec:shytaualg}).  We outline the basic theory of
hyper-Tauberian algebras, below.

\subsection{Amenability}  If $\fA$ is a Banach algebra, we let 
$\fA\what{\otimes}\fA$ denote the projective tensor product of $\fA$
with itself.  Following B. E. Johnson
\cite{johnson1}, we say $\fA$ is {\it amenable}
if it admits a {\it bounded approximate diagonal} (b.a.d.):  a bounded net
$(\mu_\alp)\subset\fA\what{\otimes}\fA$ which satisfies
\[
m(\mu_\alp)a,am(\mu_\alp)\to a\aand
a\mult\mu_\alp-\mu_\alp\mult a\to 0
\]
for $a\iin\fA$, where $m:\fA\what{\otimes}\fA\to\fA$ is the 
multiplication map, and the module actions of $\fA$
on $\fA\what{\otimes}\fA$ are given on elementary tensors
by $a\mult(b\otimes c)=(ab)\otimes c\aand (b\otimes c)\mult a
=b\otimes(ca)$.  As shown in \cite{johnson1}, amenability is
equivalent to the existence of a {\it virtual diagonal}:  an
element $M\iin(\fA\what{\otimes}\fA)^{**}$ such that
\[
a\mult M=M\mult a\aand am^{**}(M)=m^{**}(M)a=a
\]
for $a\iin\fA$, where the module actions of $\fA$ on
$(\fA\what{\otimes}\fA)^{**}$ and $\fA^{**}$,
are the second adjoints of the module actions of $\fA$
on $\fA\what{\otimes}\fA$ and $\fA$, respectively, and
$m^{**}$ is the second adjoint of the multiplication map.

We can quantify amenability via the {\it amenability} constant,
which was defined in \cite{johnson}.  Let
\[
\am{\fA}=\inf\left\{\sup_\alp\norm{\mu_\alp}:
(\mu_\alp)\text{ is a b.a.d.\ for }\fA\right\}
\]
where we allow the infemum of an empty set be $\infty$.  

The above definition is equivalent to a cohomological one:
$\fA$ is amenable if every derivation into a dual Banach $\fA$-bimodule
is inner; see \cite{johnsonM} for more on this.  We say
$\fA$ is {\it weakly amenable} if every bounded derivation into
$\fA^*$ is inner.  If $\fA$ is commutative, this is equivalent to having
every bounded derivation into any symmetric bimodule be inner;
see \cite{badecd}.  We will not directly conduct any computations
with derivations.  We note the important fact that
$\bloneg$ is amenable exactly when $G$ is an amenable group 
\cite{johnsonM}.

\subsection{The hyper-Tauberian property}\label{ssec:hypertauberian}
Let $\fA$ be a commutative semisimple Banach algebra.
Suppose $\fA$ is regular on its spectrum $\fX$;
we regard $\fA$ as an algebra of functions on $\fX$.
If $\vphi\in\fA^*$ we define
\[
\supp{\vphi}=\left\{\chi\in \fX:
\begin{matrix} \text{for every neighbourhood }U\oof \chi 
\text{ there is }f \\
\iin\fA\text{ such that }\supp{f}\subset U\aand \vphi(f)\not=0
\end{matrix}\right\}.
\]
A linear operator $T:\fA\to\fA^*$ is said to be {\it local} if 
\[
\supp{Tf}\subset\supp{f}
\]
for every $f\iin\fA$.  We say $\fA$ is {\it  hyper-Tauberian} if 
every bounded local operator $T:\fA\to\fA^*$ is an 
$\fA$-module map. This concept
was developed by the second named author \cite{samei} to study
the reflexivity of the (completely bounded) derivation space of 
$\fA$.  However, it has nice applications to weak amenability and
spectral synthesis problems, which we summarise.

\begin{hypertaubprop}\label{theo:hypertaubprop}
If $\fA$ is hyper-Tauberian then

{\bf (i)} $\fA$ is weakly amenable;

{\bf (ii)} finite subsets of $\fX$ are sets of spectral synthesis; and

{\bf (iii)} if $\fA\what{\otimes}\fA$ is semi-simple, then
$\{(\chi,\chi):\chi\in\fX\}$ is a set of local synthesis for
that algebra, and hence is a set of spectral synthesis when
$\fA$ has a bounded approximate identity.
\end{hypertaubprop}

See \cite{samei} Theorem 5, Corollary 8 and Theorem 6, for the proof.

\section{Compact groups}

\subsection{Notation}  In this section we let $G$ denote a compact group.
Let $\what{G}$ denote the set of equivalence classes
of irreducible representations of $G$.  By standard abuse of notation, we
will use $\what{G}$ to denote a set or representatives, one from each 
equivalence class.  We let $d_\pi$ denote the dimension of $\pi$.
We let $\zmeasg$ denote the centre of the measure algebra.

Let for $\pi\iin\what{G}$,
\[
\chi_\pi=\mathrm{Tr}\pi(\cdot)\aand \psi_\pi=\frac{1}{d_\pi}\chi_\pi
\]
so $\chi_\pi$ is the character of $\pi$ and $\psi_\pi$ the normalised character
with $\psi_\pi(e)=1$.  If $\mu\in\measg$ we let
\[
\hat{\mu}(\pi)=\int_G\bar{\pi}(s)d\mu(s)\in\fB(\fH_\pi).
\]
If $\mu\in\zmeasg$ then it is well-known, and strighforwrd to compute that
\[
\hat{\mu}(\pi)=\int_G \bar{\psi}_\pi d\mu\mult I_{\fH_\pi}
\]
and we then let
\begin{equation}\label{eq:zmeascoeff}
\chat{\mu}(\pi)=\int_G \bar{\psi}_\pi d\mu.
\end{equation}
We note that $f\mapsto\chat{f}$ is the Gelfand transform on $\zbloneg$.

\subsection{Some functorial properties of the centre of the group
algebra.}

We recall, as observed in \cite[Prop.\ 1.5]{mosak0}, that the map
\begin{equation}\label{eq:projection}
P=P_G:\bloneg\to\zbloneg,\quad Pf(s)=\int_Gf(tst^{-1})dt
\end{equation}
is a surjective quotient map.  

\begin{tensorprod}\label{prop:tensorprod}
$\zbloneg\what{\otimes}\zbloneg\cong\zblonegg$.
\end{tensorprod}

\proof  This follows from the fact that $P_G$ 
is a surjective quotient map and that in the identification
$\bloneg\what{\otimes}\bloneg\cong\blonegg$, we have that
$P_G\otimes P_G=P_{G\times G}$.  \endpf

If $N$ is a closed normal subgroup of $G$, we have a map
\[
T_N:\fC(G)\to\fC(G/N),\quad T_Nf(sN)=\int_Nf(sn)dn
\]
for every $sN\iin G/N$.  
This map extends to a surjective quotient map from $\bloneg$ to $\blonegn$
which we again denote $T_N$.  See \cite[Thm.\ 3.5.4]{reiters}.

\begin{quotient}\label{prop:quotient}
$T_N\bigl(\zbloneg\bigr)=\zblonegn$ and $T_N:\zbloneg\to\zblonegn$
is a surjective quotient map.
\end{quotient}

\proof 
It is sufficient to verify that
\[
T_N\comp P_G=P_{G/N}\comp T_N
\]
since each are surjective quotient maps.
For $f\in\fC(G)$ we have for $s\iin G$, using Weyl's integral formula, that
\begin{align*}
T_N\comp P_Gf(sN)&=\int_N\int_Gf(tsnt^{-1})dt\,dn \\
&=\int_N\int_{G/N}\int_Nf(tn'sn(tn')^{-1})dn'\,dtN\,dn \\
&=\int_{G/N}\int_N\int_Nf(tn'sn(tn')^{-1})dn'\,dn\,dtN \\
&=\int_{G/N} T_Nf(tst^{-1}N)dtN=P_{G/N}\comp T_Nf(sN).
\end{align*}
Since $\fC(G)$ is dense in $\bloneg$ we are done.  \endpf

\begin{quotient1}\label{cor:quotient1}
If $N$ is a closed normal subgroup of $G$ then
$\am{\zbloneg}\geq\am{\zblonegn}$.  In particular,
if $\zbloneg$ is amenable, then $\zblonegn$ is amenable.
\end{quotient1}

\proof If $(\mu_\alp)\subset\zbloneg\what{\otimes}\zbloneg$ 
is an approximate diagonal for $\zbloneg$, then
it is a standard fact that $\bigl(T_N\otimes T_N(\mu_\alp)\bigr)$
is an approximate diagonal for $\zblonegn$.  \endpf

Let $\til{m}:\measgg\to\measg$ be given by
\[
\int_Gud\til{m}(\mu)=\int_{G\times G}u(st)d\mu(s,t)\ffor u\iin\fC(G).
\]
Then $\til{m}(\mu\otimes\nu)=\mu\con\nu$ and, $\til{m}$ is the weak* continuous
extension of the multiplication map 
$m:\bloneg\what{\otimes}\bloneg\cong\blonegg\to\bloneg$.

\begin{weakstarclosed}\label{prop:weakstarclosed}
{\bf (i)} $\zmeasg$ is weak* closed in $\measg$ and $\zbloneg$
is weak* dense in $\zmeasg$.

{\bf (ii)} $\til{m}(\zmeasgg)=\zmeasg$ and $\til{m}:\zmeasgg\to\zmeasg$
is a homomorphism.
\end{weakstarclosed}

\proof {\bf (i)} The product on $\measg$ is well-known to be
weak* continuous in each variable.  Hence if $(\mu_\alp)$ is
a net contained in $\zmeasg$ with weak* limit point $\mu$, then
for each $\nu\iin\measg$ we have
\[
\nu\con\mu=\text{w*-}\lim_\alp \nu\con\mu_\alp
=\text{w*-}\lim_\alp \mu_\alp\con\nu=\nu\con\mu
\]
so $\mu\in\zmeasg$.  The map $\til{P}:\measg\to\zmeasg$ given by
\[
\int_G ud\til{P}(\mu)=\int_G \int_G u(s^{-1}ts)d\mu(t) ds
\]
for $u\in\fC(G)$, is a weak* continuous extension of the map 
$P=P_G$ in (\ref{eq:projection}).  Moreover, $\til{P}$ is a quotient map
onto $\zmeasg$, being the adjoint of the injection $\zcontg\hookrightarrow\fC(G)$,
where $\zcontg$ is the convolutive centre of $\fC(G)$.
Hence if $\mu\in\zmeasg$ and $(f_\alp)$ is a net from $\bloneg$ with
w*-$\lim_\alp f_\alp=\mu$, then w*-$\lim_\alp P(f_\alp)=P(\mu)=\mu$.

{\bf (ii)} Suppose that $\mu\in\zmeasgg$.  Then for any $u\iin\fC(G)$ and
$s\iin G$ we have
\begin{align*}
\int_G u\,d&\bigl(\del_s\con\til{m}(\mu)\con\del_{s^{-1}}\bigr)
=\int_G u(s^{-1}ts)d\til{m}(\mu)(t)=\int_Gu(s^{-1}xys)d\mu(x,y) \\
&\quad=\int_{G\times G}u(s^{-1}xss^{-1}ys)d\mu(x,y)
=\int_{G\times G}u(xy)d\mu(x,y)=\int_Gu\,d\til{m}(\mu)
\end{align*}
so $\til{m}(\mu)\in\zmeasg$.  Hence $\til{m}(\zmeasgg)\subset\zmeasg$.
Now if $\mu\in\zmeasgg$ and $\nu\in\measgg$ then for $u\iin\fC(G)$ we have
\begin{align*}
\int_Gu\,d\til{m}(\mu\con\nu)&=\int_{G\times G}u(st)d(\mu\con\nu)(s,t) \\
&= \int_{G\times G}\int_{G\times G}u(ss'tt')d\mu(s,t)d\nu(s',t') \\
&=\int_{G\times G}\int_{G\times G}u(sts't')d\mu(s,t)d\nu(s',t'), \\
&\qquad\qquad\qquad\text{ since }\del_{(e,s')}\con\mu\con\del_{(e,{s'}^{-1})}=\mu 
\text{ for any }s'\\
&=\int_G\int_G u(xy)d\til{m}(\mu)(x)d\til{m}(\nu)(y)=\int_G u\,d\bigl(\til{m}(\mu)\con
\til{m}(\nu)\bigr).
\end{align*}
Observe that we actually proved that
$\til{m}$ is a (left) $\zmeasgg$-module map.
Since $\til{m}(\mu\otimes\del_e)=\mu$, it follows that
$\til{m}(\zmeasgg)=\zmeasg$.
\endpf

\subsection{Approximate diagonals for centres of compact group algebras}
We note that $\ztrig=\spn\{\chi_\pi:\pi\in\what{G}\}$ is dense in
$\zbloneg$.  To see this, we first recall that the set $\trig=
\{\pi_{ij}:i,j=1,\dots,d_\pi,\pi\in\what{G}\}$ of matrix coefficients
is dense in $\bloneg$.  It is easily checked that 
$P\pi_{ij}=\psi_\pi=d_\pi^{-1}\chi_\pi$ for each $\pi_{ij}$
where $P$ is the map defined in (\ref{eq:projection}).
Then if $(u_n)\subset\trig$ converges to $f\iin\zbloneg$,
we have $\lim_nPu_n=Pf=f$.

\begin{nicebai}\label{lem:nicebai}
There exists a net $(f_\beta)$ in $\ztrig$ such that
$(f_\beta)$ is a bounded approximate identity for $\bloneg$.
Moreover, if for each $\beta$ we have
\[
f_\beta=\sum_{\pi\in\what{G}}a^\beta_\pi\chi_\pi
\]
where $a^\beta_\pi=0$ except for finitely many elements $\pi$, then
for each $\pi\iin\what{G}$ we have
\[
\lim_{\beta}a^\beta_\pi=d_\pi.
\]
\end{nicebai}

\proof Let $(U)$ be a base of neighbourhoods of the identity
in $G$, each  invariant for inner automorphisms.  Then
$(e_U)=\left(\frac{1}{\lam(U)}1_U\right)$ is a central approximate
identity for $\bloneg$.   Since $\ztrig$ is dense in $\zbloneg$
we can find for each $\eps>0$ and $U$ as above, $f_{\eps,U}\in\ztrig$
such that $\norm{f_{\eps,U}-e_U}_1<\eps$.  Then $(f_\beta)=(f_{\eps,U})$
is the desired bounded approximate identity.

Since for each $\pi\iin\what{G}$ we have
\[
\frac{a^\beta_\pi}{d_\pi}\chi_\pi=f_\beta\con\chi_\pi
\overset{\beta}{\longrightarrow}\chi_\pi
\]
it follows that $\lim_\beta a^\beta_\pi=d_\pi$.  \endpf

We recall from \cite[(27.43)]{hewittrII} that $\what{G\cross G}
=\{\pi\cross\sig:\pi,\sig\in\what{G}\}$.

\begin{diagonal}\label{theo:diagonal}
Let $G$ be a compact group and $(f_\beta)$ be as in Lemma \ref{lem:nicebai},
above.  For each $\beta$ define
\[
\mu_\beta=\sum_{\pi\in\what{G}}(a_\pi^\beta)^2\chi_\pi\otimes\chi_\pi
\iin \zbloneg\what{\otimes}\zbloneg.
\]
Then $(\mu_\beta)$ is an approximate diagonal for $\zbloneg$.
Moreover, the following are equivalent

{\bf (i)} $(\mu_\beta)$ is bounded;

{\bf (ii)} $\zbloneg$ is amenable; and

{\bf (iii)} there is a measure $\mu\iin\zmeasgg$ which satisfies
\begin{equation}\label{eq:amenzmeas}
\chat{\mu}(\pi\cross\sig)=\del_{\pi,\sig}
\end{equation}
where $\delta$, in this context, is the Kronecker symbol.   For such $\mu$ we
have that $\til{m}(\mu)=\del_e$ and $(f\otimes\del_e)\con \mu=
\mu\con(\del_e\otimes f)$ for $f\iin\zbloneg$.
\end{diagonal}

Note that we thus have that $\zbloneg$ is {\it psuedo-amenable},
in the sense defined in \cite{ghahramaniz}.

\medskip
\proof It is clear that for each $\pi\iin\what{G}$ and each $\beta$ we have
$\chi_\pi\mult\mu_\beta=\mu_\beta\mult\chi_\pi$.  Since
$\ztrig$ is dense in $\zbloneg$, it follows that
$f\mult \mu_\beta=\mu_\beta\mult f$ for each $f\iin\zbloneg$ too.  Also
\[
m(\mu_\beta)=\sum_{\pi\in\what{G}}\frac{(a_\pi^\beta)^2}{d_\pi}\chi_\pi
=f_\beta\con f_\beta
\]
so $\bigl(m(\mu_\beta)\bigr)$ is a bounded approximate identity.
Thus $(\mu_\beta)$ is an approximate diagonal for $\zbloneg$.
It is immediate that {\bf (i)} $\implies$ {\bf (ii)}.

{\bf (ii)} $\implies$ {\bf (i)}.
If we suppose that $\zbloneg$ is amenable, it admits a bounded approximate
diagonal $(\mu'_\gam)$.  We may assume $(\mu'_\gam)$ is weakly Cauchy, i.e.\ 
it converges to a virtual diagonal $M$ in 
$\bigl(\zbloneg\what{\otimes}\zbloneg\bigr)^{**}$.
With $(f_\beta)$ as in the lemma above, let 
for each $\beta$, $F_\beta=
\{\pi\in\what{G}:a^\beta_\pi\not=0\}$ and $\fA_\beta=\spn\{\chi_\pi:\pi\in F_\beta\}$.
Then $\fA_\beta\otimes\fA_\beta$ is a finite dimensional ideal in 
$\zbloneg\what{\otimes}\zbloneg$ which contains $f_\beta\otimes f_\beta$.
Then $\bigl((f_\beta\otimes f_\beta)\con \mu'_\gam\bigr)_\gam$ is a bounded net
in $\fA_\beta\otimes\fA_\beta$ with limit point $\mu'_\beta$. 
Write
\[
\mu'_\beta=\sum_{\pi,\sig\in F_\beta}c^\beta_{\pi,\sig}\chi_\pi\otimes\chi_\sig.
\]
Then for any $\pi\in F_\beta$, using that $f\con\chi_\pi=\chi_\pi\con f$, we have
\[
\chi_\pi\mult \mu'_\beta
=(f_\beta\otimes f_\beta)\con\left[\lim_\gam \chi_\pi\mult \mu'_\gam\right]
=(f_\beta\otimes f_\beta)\con\left[\lim_\gam \mu'_\gam\mult \chi_\pi\right]
=\mu'_\beta\mult\chi_\pi
\]
and thus
\[
\sum_{\sig\in F_\beta}\frac{c^\beta_{\pi,\sig}}{d_\pi}\chi_\pi\otimes\chi_\sig
=\sum_{\sig\in F_\beta}\frac{c^\beta_{\sig,\pi}}{d_\pi}\chi_\sig\otimes\chi_\pi.
\]
It follows from the orthogonality relations of the characters that 
$c^\beta_{\pi,\sig}=0$ if $\sig\not=\pi$ and hence 
\[
\mu'_\beta=\sum_{\pi\in F_\beta}\frac{c^\beta_{\pi,\pi}}{d_\pi}\chi_\pi\otimes\chi_\pi.
\]
Since $m(\mu'_\beta)=m(f_\beta\otimes f_\beta)
\con\lim_\gam m(\mu'_\gam)=f_\beta\con f_\beta$
we obtain
\[
\sum_{\pi\in F_\beta}\frac{c^\beta_{\pi,\pi}}{d_\pi}\chi_\pi
=\sum_{\pi\in F_\beta}\frac{(a^\beta_\pi)^2}{d_\pi}\chi_\pi
\]
and thus $c^\beta_{\pi,\pi}=(a^\beta_\pi)^2$ for each $\pi\iin F_\beta$.
Then for each $\beta$ we have $\mu_\beta=\mu'_\beta$, so
\[
\norm{\mu_\beta}=\norm{\mu'_\beta}\leq\norm{f_\beta}_1^2
\sup_\gam\norm{\mu'_\gam}
\]
and, since $(f_\beta)$ is bounded, $(\mu_\beta)$ is bounded too.

{\bf (i)} $\implies$ {\bf (iii)}.
Using Proposition \ref{prop:tensorprod} we identify $(\mu_\beta)$
as a bounded net $\measgg$.  It thus has a weak* cluster point
$\mu$.  We note that $\mu$ is, in fact, a limit point.  Indeed,
$\trigg$ is uniformly dense in $\fC(G\cross G)$, and if $u\in\trigg$
it is clear that
\begin{align*}
\sum_{\pi\in\what{G}}&d_\pi^2\int_{G\times G} u(s,t)\chi_\pi(s)\chi_\pi(t)d(s,t) \\
&=\lim_{\beta}\sum_{\pi\in\what{G}}(a^\beta_\pi)^2
\int_{G\times G} u(s,t)\chi_\pi(s)\chi_\pi(t)d(s,t) 
=\lim_\beta\dpair{u}{\mu_\beta}
\end{align*}
as all sums in the expression are finite.  Moreover,
the above expression must be $\int_{G\times G}u(s,t)d\mu(s,t)$.
By Proposition \ref{prop:weakstarclosed}, $\mu\in\zmeasgg$.
By (\ref{eq:zmeascoeff}) we find
\[
\chat{\mu}(\pi\cross\sig)
=\frac{1}{d_\pi d_\sig}\int_{G\times G}\wbar{\chi_\pi(s)\chi_\sig(t)}d\mu(s,t)
=\del_{\pi,\sig}.
\]

Let $R:\fC(\what{G}\cross\what{G})\to\fC(\what{G})$ be the map of restriction
to the diagonal:  $Ru(\pi)=u(\pi,\pi)$.  Note that for any $\nu\iin\zmeasgg$,
$\bigr(\til{m}(\nu)\bigl)^\smallfrown=R\chat{\nu}$.  Thus we have for 
any $\pi\iin\what{G}$
\[
\bigl(\til{m}(\mu)\bigr)^\smallfrown(\pi)=R\chat{\mu}(\pi)=1=\chat{\del_e}(\pi)
\]
so $\til{m}(\mu)=\del_e$.  Also, if $f\in\zbloneg$ and $\pi,\sig\in\what{G}$ then,
$(f\otimes\del_e)^\smallfrown(\pi\cross\sig)=\chat{f}(\pi)$ while
$(\del_e\otimes f)^\smallfrown(\pi\cross\sig)=\chat{f}(\sig)$.  It follows that
\[
(f\otimes\del_e)^\smallfrown\chat{\mu}(\pi\cross\sig)=
f(\pi)\del_{\pi,\sig}=\chat{\mu}(\del_e\otimes f)^\smallfrown(\pi\cross\sig)
\]
so it follows that $(f\otimes\del_e)\con\mu=\mu\con(\del_e\otimes f)$.

{\bf (iii)} $\implies$ {\bf (ii)}.  Let $(f_\alp)$ be any bounded approximate
identity in $\zblonegg$.  We will show that any weak* cluster point $M$ of
$(\mu\con f_\alp)$ in $\zblonegg^{**}$ is a virtual diagonal.
We may assume $M$ is a limit point.  First, if $f\in\zbloneg$ we have
\[
f\mult M =\lim_\alp (f\otimes\del_e)\con\mu\con f_\alp=
\lim_\alp\mu\con(\del_e\otimes f)\con f_\alp 
= \lim_\alp\mu\con f_\alp\con (\del_e\otimes f)=M\mult f.
\]
Second, we note it follows from Proposition \ref{prop:weakstarclosed} that
$\bigl(m(f_\alp)\bigr)$ is a bounded approximate identity for $\zbloneg$.
We let $E$ be any weak* cluster point of $\bigl(m(f_\alp)\bigr)$, which we may
consider to be a limit point.  We then have, again by Proposition 
\ref{prop:weakstarclosed}, and using $\til{m}(\mu)=\del_e$, that
\[
m^{**}(M)=\lim_\alp m(\mu\con f_\alp)=\lim_\alp \til{m}(\mu)\con m(f_\alp)
=\lim_\alp m(f_\alp)=E.
\]
It is clear that $f\mult E=E\mult f=f$ for $f\in\zbloneg$.  Thus $M$ is a virtual diagonal.
\endpf

Note that if $G$ is abelian, then $\mu$ is the Haar measure of the
anti-diagonal subgroup $A=\{(s,s^{-1}):s\in G\}$.  Indeed, if we denote the
latter by $\lam_A$ then we have for $\chi,\psi\iin\what{G}$
\[
\chat{\lam}_A(\chi\cross\psi)=\int_G \wbar{\chi(s)\psi(s^{-1})}ds
=\int_G \wbar{\chi(s)}\psi(s)ds=\del_{\chi,\psi}=\chat{\mu}(\chi\cross\psi)
\]
and hence $\mu=\lam_A$.  Though the definition of $\lam_A$, as above,
makes sense for any compact group, it forms a central measure only 
when $G$ is abelian.

Suppose $d_G=\sup_{\pi\in\what{G}}d_\pi<\infty$.  Then for $u,v\iin\trig$
we use the Cauchy-Schwarz inequality 
and Bessel's inequlity on the ortho-normal set $\{\chi_\pi:\pi\in\what{G}\}$
to see that for the approximite diagonal  $(\mu_\beta)$ in the theorem above
we have
\begin{align*}
\left|\lim_\beta\int_{G\times G}u(s)v(t)\mu_\beta(s,t)dsdt\right|
&=\left|\sum_{\pi\in\what{G}}d_\pi^2\int_Gu(s)\chi_\pi(s)ds\cdot\int_Gv(t)\chi_\pi(t)dt 
\right| \\
&\leq
d_G^2\sum_{\pi\in\what{G}}|\inprod{u}{\chi_{\bar{\pi}}}||\inprod{v}{\chi_{\bar{\pi}}}|
\\
&\leq d_G^2\norm{u}_2\norm{v}_2\leq d_G^2\norm{u}_\infty\norm{v}_\infty.
\end{align*}
Since $\trig$ is dense in $\fC(G)$, it follows that $(\mu_\beta)$
converges to a {\it bimeasure} in the terminology of \cite{grahams}, 
i.e.\ an element $\mu$ of 
$\bigl(\fC(G)\what{\otimes}\fC(G)\bigr)^*$.  
Conjecture \ref{conj:mainconj}, if true, would further imply that
$\mu\in\measgg$.

\subsection{Connected groups}\label{ssec:conngroup}

\begin{conngroup}\label{theo:conngroup}
If $G$ is a non-abelian connected compact group, then $\zbloneg$
is not amenable.
\end{conngroup}

\proof There is a family $\{G_i\}_{i\in I}$ 
of compact connected Lie groups, at least one of which is simple (in the sense 
of Lie groups) with finite centre, such that
\[
G\cong \left(\prod_{i\in I}G_i\right)/A
\]
where $A$ is a central subgroup of $P=\prod_{i\in I}G_i$.
See \cite[6.5.6]{price}, for example. Hence $G$ admits, as a quotient
\[
\prod_{i\in I}G_i/\z{G_i}\cong
P/\z{P}\cong(P/A)/(\z{P}/A).
\]
Let $i_0$ be so $G_{i_0}$ is simple with finite centre.  Then 
$G_{i_0}/\z{G_{i_0}}$ is simple with trivial centre.  Hence
there is a closed normal subgroup $N$ of $G$ such that
$G/N$ is a simple Lie group with trivial centre.  By \cite[Lem.\ 9.1]{rider}
we obtain ``Condition I" on $G/N$, which is the property that
\[
\lim_{d_\pi\to\infty}\psi_\pi(sN)= 0\ffor sN\in G/N\setdif\{eN\}.
\]
Hence, there is a sequence $\{\pi_n\}_{n=1}^\infty\subset\what{G}$
such that
\begin{equation}\label{eq:weakrider}
\psi_n(s)=1\ffor s\iin N\aand \lim_{n\to\infty}\psi_n(s)=0\ffor s\in G\setdif N
\end{equation}
where $\psi_n=\psi_{\pi_n}$.  Indeed, choose any sequence of representations
$\{\til{\pi}_n\}_{n=1}^\infty\subset\what{G/N}$ where $\lim_{n\to\infty}d_{\til{\pi}_n}
=\infty$, and let $\pi_n=\til{\pi}_n\comp q$ where $q:G\to G/N$ is the quotient map.

If it were the case that $\zbloneg$ were amenable, then we would obtain
$\mu\iin\zmeasgg$ as in (\ref{eq:amenzmeas}).  Let us see that
the existence of such $\mu$ gives a contradiction.
Let  $N$ and $(\psi_n)$ be as in (\ref{eq:weakrider}).  Define
two sequences $(u_n)$ and $(v_n)$ of functions on $G\cross G$ by
\[
u_n=\psi_n\otimes\psi_n\aand v_n=\psi_n\otimes\psi_{n+1}.
\]
Then $(u_n)$ and $(v_n)$ are bounded sequences with
\[
\lim_{n\to\infty}u_n(s,t)=\lim_{n\to\infty}v_n(s,t)=\begin{cases}
1 &\iif (s,t)\in N\cross N \\ 0 & \iif (s,t)\not\in N\cross N.\end{cases}
\]
Hence it follows from the Lebesgue dominted convergence theorem that
\begin{equation}\label{eq:ldct}
\lim_{n\to\infty} \int_{G\times G}u_nd\mu=\mu(N\cross N)=
\lim_{n\to\infty} \int_{G\times G}v_nd\mu.
\end{equation}
However, by (\ref{eq:amenzmeas}) we have that
\[
\int_{G\times G}u_nd\mu=\chat{\mu}(\bar{\pi}_n\cross\bar{\pi}_n)=1
\;\text{ while }\;
\int_{G\times G}v_nd\mu=\chat{\mu}(\bar{\pi}_n\cross\bar{\pi}_{n+1})=0
\]
which contradicts (\ref{eq:ldct}).  \endpf

\subsection{Products of finite groups}\label{ssec:prodgroup}
Let $G$ be a finite group.  We will treat $G$ as a compact group so we have
normalised Haar integral: $\int_G f=\frac{1}{|G|}\sum_{s\in G}f(s)$.
Then it is well known that
\begin{equation}\label{eq:zblonechar}
\zbloneg=\spn\{\chi_\pi:\pi\in\what{G}\}
\end{equation}
Moreover, if we let for any $x\iin G$, $C_x=\{sxs^{-1}:s\in G\}$ denote
the conjugacy class, and $\conjg=\{C_x:x\in G\}$, then since elements
of $\zbloneg$ are constant on conjugacy classes we have
\begin{equation}\label{eq:zbloneconj}
\zbloneg=\spn\{1_C:C\in\conjg\}
\end{equation}
where $1_C$ is the indicator function of $C$.  We will let
$f(C)=f(x)$ where $C=C_x$, for $f\in\zbloneg$.

\begin{fingroup}\label{theo:fingroup}
If $G$ is a finite group, then $\zbloneg$ has unique diagonal
and we have 
\[
\am{\zbloneg}=\frac{1}{|G|^2}\sum_{C,C'\in\conjg}
\left|\sum_{\pi\in\what{G}}d_\pi^2\wbar{\chi_\pi(C)}\chi_\pi(C')\right||C||C'|.
\]
\end{fingroup}

\proof That
\[
\mu=\sum_{\pi\in \what{G}}d_\pi^2\chi_\pi\otimes\chi_\pi
\]
is the unique diagonal for $\zbloneg$ follows from the proof of
Theorem \ref{theo:diagonal}.  However, using the relations
$\chi_\pi\con\chi_\sig=\del_{\pi,\sig}d_\pi^{-1}\chi_\pi$ for
$\pi,\sig\iin\what{G}$, that $\mu$ is a diagonal is easily verified manually 
using (\ref{eq:zblonechar}).  The uniqueness of the diagonal in any 
amenable finite dimensonal commutative algebra has been observed
in \cite[Prop.\ 0.2]{ghandeharihs}.

If $C\in\conjg$ with $C=C_x$, we let $\wbar{C}=C_{x^{-1}}$.  The operation
$C\mapsto\wbar{C}$ is an involution on $\conjg$.  If $\pi\in\what{G}$
then $\chi_\pi(\wbar{C})=\wbar{\chi_\pi(C)}$.
We appeal to (\ref{eq:zbloneconj}) to obtain
\begin{align*}
\mu&=\sum_{\pi\in \what{G}}d_\pi^2
\left(\sum_{C\in\conjg}\chi_\pi(C)1_C\right)\otimes
\left(\sum_{C'\in\conjg}\chi_\pi(C')1_{C'}\right) \\
&=\sum_{\pi\in \what{G}}d_\pi^2
\left(\sum_{C\in\conjg}\wbar{\chi_\pi(C)}1_C\right)\otimes
\left(\sum_{C'\in\conjg}\chi_\pi(C')1_{C'}\right) \\
&=\sum_{C,C'\in\conjg}
\left(\sum_{\pi\in\what{G}}d_\pi^2\wbar{\chi_\pi(C)}\chi_\pi(C')\right)
1_C\otimes 1_{C'}.
\end{align*}
We then compute $\am{\zbloneg}=\norm{\mu}_1$ to finish.  \endpf

\begin{fingroup1}\label{cor:fingroup1}
If $G$ is a non-abelian finite group, then $\am{\zbloneg}>1$.
\end{fingroup1}

\proof Letting $C'=C$ we obtain lower bound
\begin{align*}
\am{\zbloneg}&\geq \frac{1}{|G|^2}\sum_{C\in\conjg}\sum_{\pi\in\what{G}}d_\pi^2
|\chi_\pi(C)|^2|C|^2 \\
&=\frac{1}{|G|}\sum_{\pi\in\what{G}}d_\pi^2\sum_{C\in\conjg}|C||\chi_\pi(C)|^2
\frac{|C|}{|G|}.
\end{align*}
Since $G$ is nonabelian we have $|C|>1$ for some conjugacy class $C$.
Moreover, there is some $\pi$ so $\chi_\pi(C)\not=0$.  Thus we find
\begin{align*}
\am{\zbloneg}&>
\frac{1}{|G|}\sum_{\pi\in\what{G}}d_\pi^2\sum_{C\in\conjg}|\chi_\pi(C)|^2
\frac{|C|}{|G|} \\
&= \frac{1}{|G|}\sum_{\pi\in\what{G}}d_\pi^2\norm{\chi_\pi}_2^2=1
\end{align*}
since $\norm{\chi_\pi}_2=1$ and $\sum_{\pi\in\what{G}}d_\pi^2=|G|$.  \endpf

Let us take a second look at the proof of the above corollary.
The Schur orthogonality relations tell us that the $\what{G}\cross\conjg$
matrix 
\[
U=\left[\frac{|C|^{1/2}}{|G|^{1/2}}\chi_\pi(C)\right]
\]
is unitary.  
Letting $C'=C$ we obtain lower bound
\begin{align*}
\am{\zbloneg}&\geq \frac{1}{|G|}\sum_{C\in\conjg}\sum_{\pi\in\what{G}}d_\pi^2
|\chi_\pi(C)|^2\frac{|C|}{|G|}|C| \\
&=\frac{1}{|G|}\norm{\mathrm{diag}(d_\pi)U\mathrm{diag}(|C|^{1/2})}_2^2
\end{align*}
where $\norm{\cdot}_2$ denotes the Hilbert-Schmidt norm.  
 {\it Is it possible to get a lower 
estimate in terms of $\max_{\pi\in\what{G}}d_\pi$?} 
If so, Conjecture \ref{conj:mainconj} may be shown to hold for compact totally 
disconnected groups which do not admit an open abelian subgroup.

\begin{prodfingroup}\label{theo:prodfingroup}
If $G=\prod_{i=1}^\infty G_i$ where each $G_i$ is a nonabelian finite group,
then $\zbloneg$ is not amenable.
\end{prodfingroup}

\proof  For each $i$, the diagonal $\mu_i$  for $\zblone{G_i}$ promised by
Theorem \ref{theo:fingroup} is an idempotent in $\zblonegg$.  Hence
by \cite[Thm.\ 5.3]{rider}, there is a constant $\del>0$ -- in fact $\del\geq 1/700$ --
for which
\[
\am{\zblone{G_i}}\geq 1+\del
\]
for each $i$.  Since $G$ admits, for each $n$, $G_{(n)}=\prod_{i=1}^nG_i$ as a 
quotient, we have that
\[
\am{\zbloneg}\geq\am{\zblone{G_{(n)}}}=\prod_{i=1}^n\am{\zblone{G_i}}\geq
(1+\del)^n.
\]
Hence we have that $\am{\zbloneg}=\infty$ and $\zbloneg$ is not amenable.
\endpf

\subsection{An amenable example}\label{ssec:amenex}
The following example further illustrates Conjecture \ref{conj:mainconj}.

Let $G=\Tee\ltimes\Zee_2$ where $\Tee=\{s\in\Cee:|s|=1\}$ and $\Zee_2=\{-1,1\}$.
The group law and inverse are given by
\[
(s,a)(t,b)=(st^a,ab)\aand(s,a)^{-1}=(s^{-a},a)
\]
for $(s,a),(t,b)\iin G$.  An application of the ``Mackey machine", see 
\cite[Sec.\ 6.6]{folland} for example, gives us $\what{G}=\{1,\sig,\pi_n:n\in\En\}$
where
\[
1(s,a)=1,\;\sig(s,a)=a\aand
\pi_n(s,a)=\begin{bmatrix} s^n & 0 \\ 0 & s^{-n}\end{bmatrix}
\begin{bmatrix} 0 & 1\\ 1 & 0\end{bmatrix}^{(1-a)/2}
\]
for $(s,a)\iin G$. It follows that we have normalised characters $1,\sig$ and
\[
\psi_{\pi_n}(s,a)=\begin{cases} \frac{1}{2}(s^n+s^{-n}) &\iif a=1 \\
\;\; \quad\;\; 0&\iif a=-1\end{cases}
\]
for $(s,a)\iin G$.  We note that all of the calculations thus far, and hence the
next proposition, also hold if $\Tee$ is replaced by any compact abelian group
$T$ admiting only $1$ as a real character, i.e.\ for $\chi\iin\what{T}$,
$\wbar{\chi}=\chi$ implies $\chi=1$.  For sake of concreteness, we will
continue with $T=\Tee$.

\begin{finabexample}\label{prop:finabexample}
For $G=\Tee\ltimes\Zee_2$, $\zbloneg$ is amenable.
\end{finabexample}

\proof Let $\mu=1\otimes 1+\sig\otimes\sig-2(1+\sig)\otimes(1+\sig)+\nu$, 
where $\nu=\lam_D+\lam_A$, the sum of the Haar measures on the subgroups
of $G\cross G$ given by $D=\left\{\bigl((t,1),(t,1)\bigr):t\in\Tee\right\}$ and
$A=\left\{\bigl((t,1),(t^{-1},1)\bigr):t\in\Tee\right\}$, each normalised
to have total mass $1$.  We note that
$\nu\in\zmeasg$ since for $(s,a)\iin G$ we have
\begin{align*}
\del_{((s,a),(t,b))}&\con\nu\con\del_{((s^{-a},a),(t^{-b},b))} \\
=&\del_{((1,a),(1,b))}\con\del_{((s^a,1),(t^b,1))}\con\nu\con
\del_{((s^{-a},1),(t^{-b},1))}\con\del_{((1,a),(1,b))}
=\nu.
\end{align*}
Thus $\mu\in\zmeasg$.  Now for $\pi,\rho\iin\what{G}$ we have
\begin{align*}
\chat{\mu}(\pi\cross\rho)
=&\int_{G\times G}\bigl(1\otimes 1+\sig\otimes\sig-2(1+\sig)\otimes(1+\sig)\bigr)
\mult\psi_\pi\otimes\psi_\rho \\
&+\int_\Tee 
\bigl(\psi_\pi(s,1)\psi_\rho(s,1)+\psi_\pi(s,1)\psi_\rho(s^{-1},1)\bigr)ds \\
&=(1)+(2)
\end{align*}
where
\begin{align*}
(1)&=\begin{cases} -1 &\iif (\pi,\rho)=(1,1),(\sig,\sig) \\
-2 &\iif (\pi,\rho)=(1,\sig),(\sig,1) \\
\phantom{-}0 &\iif (\pi,\rho)=(1,\pi_n),(\pi_n,1),(\sig,\pi_n),(\pi_n,\sig) \\
\phantom{-0} &\phantom{\iif (\pi,\rho)= }
\; (\pi_n,\pi_m),\, n,m\in\En\end{cases}\qquad\aand \\
(2)&=2\int_\Tee\psi_\pi(s,1)\psi_\rho(s,1)ds
=\begin{cases}2 &\iif (\pi,\rho)=(1,1),(\sig,\sig),(1,\sig),(\sig,1) \\
1&\iif (\pi,\rho)=(\pi_n,\pi_n),\,n\in\En \\
0 &\iif (\pi,\rho)=(1,\pi_n),(\pi_n,1),(\sig,\pi_n),(\pi_n,\sig) \\
\phantom{0} & \phantom{\iif (\pi,\rho)=}\; (\pi_n,\pi_m),\, n\not=m, \,n,m\in\En 
\end{cases}
\end{align*}
Thus it follows that $\mu$ satisfies (\ref{eq:amenzmeas}).  

We remark that the measure $\mu$ corresponds to the (formal) Fourier series
\[
1\otimes 1+\sig\otimes\sig-2(1+\sig)\otimes(1+\sig)
+4\left[\frac{1}{2}(1+\sig)\otimes(1+\sig)+\sum_{n=1}^\infty
\chi_{\pi_n}\otimes\chi_{\pi_n}\right]
\]
as suggested by Theorem \ref{theo:diagonal}.
The coefficient $4$, in the second term, is $1/\lam(G_{(1,1)}\cross G_{(1,1)})$, 
where $G_{(1,1)}\cong\Tee$ is the connected component of the identity in $G$.
The second term corresponds to the Fourier series for $\lam_A+\lam_D$
on $\Tee\cross\Tee$, as may be revealed by a simple computation which
we leave to the reader. 
\endpf

Let us make a few observations about $G=\Tee\ltimes\Zee_2$.
First we compute, for $s\iin\Tee\aand(t,b)\iin G$
\[
(t,b)(s,1)(t^{-b},b)=(s^b,1)\aand (t,b)(s,-1)(t^{-b},b)=(t^2s^b,-1).
\]
Hence we deduce that
\[
\conjg=\bigl\{\{(1,1)\},\{(-1,1)\},\{(s,1),(s^{-1},1)\}_{\mathrm{Im}s>0},
G_{(1,-1)}\bigr\}
\]
where $G_{(1,-1)}$ is the connected component of $(1,-1)$.  Moreover
we compute commutators
\[
[(t,b),(s,a)]=(t,b)(s,a)(t^{-b},b)(s^{-a},a)=(t^{1-a}s^{1-b},1).
\]
Letting $a=1,b=-1$ and $s,t$ be arbitrary in $\Tee$ we find, in the notation of
Section \ref{ssec:fincommutator} that $G'=[G,G_0]=G_{(1,1)}$. 
In particular, notice that the assumptions
of Theorem \ref{theo:finitederived}, below, are not necessary for $\zbloneg$
to be amenable.

Let us close by noting the following decomposition
\[
\zblone{\Tee\ltimes\Zee_2}=\mathrm{Z}_{\Zee_2}\blone{\Tee}
\oplus\Cee (1-\sig)
\]
where $\mathrm{Z}_{\Zee_2}\blone{\Tee}=\{f\in\blone{\Tee}:\check{f}
=f\}$, $\check{f}(s)=f(s^{-1})$.  We note that both of the components of this 
decomposition are closed subalgebras, but neither is an ideal.
Hence it is not apparent that $\mathrm{Z}_{\Zee_2}\mathrm{L}^1(\Tee)$
is amenable.  We show this fact in the next section.

\subsection{The hypergroup approach}\label{ssec:hypergroup}
We indicate, by way of two examples, how the problem of amenability for
$\zbloneg$ can be treated by using hypergroups.  We refer to
\cite{bloomh} for the definintion of a hypergroup $K$ and its left Haar measure
$\lam_K$, or to \cite{jewett}, where a hypergroup is referred to as a ``convos".
If $G$ is a compact group, then $K=\conjg$ is a hypergroup
\cite[8.4]{jewett}.  Since $K$ is compact and commutative, it admits a Haar measure.
Moreover we have $\zbloneg\cong\blonek$, where $\blonek$ is the hypergroup
algebra.
Such $K$ is a {\it strong hypergroup} in the sense that its character set
$\what{K}$ is a (discrete) hypergroup under pointwise multiplication.  In fact
$\what{K}$ identifies naturally with $\{\psi_\pi\}_{\pi\in\what{G}}$.

We first consider $G=\sutwo$.  By \cite[15.4]{jewett}, $\conjsutwo$ identifies
naturally with a hypergroup whose underlying set is $K=[-1,1]$.  We
will not explicitly need the convolution formula on $K$, but we will require
the formula for the Haar measure
\begin{equation}\label{eq:khaarmeas}
\int_K fd\lam_K=\frac{2}{\pi}\int_{-1}^1f(x)(1-x^2)^{1/2}dx=
\frac{2}{\pi}\int_0^\pi f(\cos\theta)\sin^2\theta d\theta
\end{equation}
where $dx,d\theta$ each denote integration with respect to Lebesgue measure,
and the (non-normalised) characters are given by
\[
\{\chi_k\}_{k=0}^\infty\qquad\wwhere\qquad 
\chi_k(\cos\theta)=\frac{\sin(k+1)\theta}{\sin\theta}.
\]
Note that $\chi_k$ is, up to identification, the character
of the unique representation of $\sutwo$ of dimension $k+1$.

\begin{sutwona}\label{theo:sutwona}
$\zblone{\sutwo}$ is not amenable.
\end{sutwona}

\proof  We first note that by Lemma \ref{lem:nicebai}, there is a bounded
approximate identity for $\blonek$, $(e_\alp)\subset\mathrm{Trig}(K)
=\spn\{\psi_k:k\in\En_0\}$ where $\En_0=\{0\}\cup\En$.  This bounded
approximate identity may be taken to be a sequence, $(e_n)$, 
and we have for each $n$, 
$e_n=\sum_{k=0}^\infty a_k^{(n)}\chi_k$
where $a_k^{(n)}=0$ for all but finitely many indices $k$.
We obtain, again from
Lemma \ref{lem:nicebai}, that
\begin{equation}\label{eq:somelimit}
\lim_{n\to\infty}a_k^{(n)}=k+1
\end{equation}
Now let $\mu_n=\sum_{k=0}^\infty(a_k^{(n)})^2\chi_k\otimes\chi_k$, so
$(\mu_n)_{n=1}^\infty$ is the approximate identity from Theorem
\ref{theo:diagonal}, and we are done once we establish $(\mu_n)$ is
not bounded.  
The using the fact that $\blonek\what{\otimes}\blonek\cong
\blone{K\cross K}$ and (\ref{eq:khaarmeas}) we have
\begin{align*}
\left(\frac{\pi}{2}\right)^2&\norm{\mu_n}_1
=\int_0^\pi\int_0^\pi
|\mu_n(\cos\theta,\cos\theta')|\sin^2\theta\sin^2\theta'd\theta d\theta' \\
&\geq\int_0^{\pi/2}\int_0^{\pi/2}\left|\sum_{k=0}^\infty
(a_k^{(n)})^2
\frac{\sin(k+1)\theta}{\sin\theta}\frac{\sin(k+1)\theta'}{\sin\theta'}
\right|\sin^2\theta\sin^2\theta'd\theta d\theta' \\
&\geq\left|\int_0^{\pi/2}\int_0^{\pi/2}
\sum_{k=0}^\infty(a_k^{(n)})^2
\sin(k+1)\theta\sin(k+1)\theta'
\sin\theta\sin\theta'd\theta d\theta' \right| \\
&=\sum_{k=0}^\infty(a_k^{(n)})^2
\left(\int_0^{\pi/2}\sin(k+1)\theta\sin\theta d\theta\right)^2 \\
&=\sum_{k=0}^\infty(a_k^{(n)})^2
\left(\frac{1}{2}\int_0^{\pi/2}\bigl(\cos k\theta-\cos(k+2)\theta\bigr)d\theta\right)^2 \\
&=\sum_{k=0}^\infty
\left(\frac{a_k^{(n)}(k+1)}{k(k+2)}\sin\frac{k\pi}{2}\right)^2 
=\sum_{j=0}^\infty
\left(\frac{a_{2j+1}^{(n)}(2j+2)}{(2j+1)(2j+3)}\right)^2
\end{align*}
Let $f_n=\left(\frac{a_{2j+1}^{(n)}(2j+2)}{(2j+1)(2j+3)}\right)_{j=0}^\infty$.
If $(\mu_n)$ is bounded, then $(f_n)$ is bounded in $\ell^2(\En_0)$, in which
cases the latter sequence has a cluster point $f$.  We have, by
(\ref{eq:somelimit}), that $f(j)=\frac{(2j+2)^2}{(2j+1)(2j+3)}$, which means that
$f$ cannot be an element of $\ell^2(\En_0)$.  Thus $(\mu_n)$ must not be bounded.
\endpf

Let us now turn out attention to $\mathrm{Z}_{\Zee_2}\blone{\Tee}$, from
the last section.  We let
\[
\psi_0=\frac{1}{2}(1+\sig)\aand \psi_n=\psi_{\pi_n}\ffor n\in\En.
\]
Then the family of all $\Zee_2$-invariant characters of 
$\mathrm{Z}_{\Zee_2}\blone{\Tee}$ is 
$\fX_{\Zee_2}(\Tee)=\{\psi_n\}_{n\in\En_0}$.  
Observe, under pointwise multiplication, that $\fX_{\Zee_2}(\Tee)$
satisfies the same multiplication rules as the cosine functions
$\{\cos(m\cdot)\}_{m\in\En_0}$, and hence
is isomorphic to the Chebychev polynomial hypergroup of the first kind
\cite{bloomh}.

There is a commuative
hypergroup $K=[-1,1]$, which is isomorphic to the double conjugacy
class hypergroup $\Tee/\!/\Zee_2$, such that $\mathrm{Z}_{\Zee_2}\blone{\Tee}
\cong\blonek$.  The Haar measure on $K$ is given by
\[
\int_K fd\lam_K=\frac{1}{\pi}\int_{-1}^1\frac{f(x)}{\sqrt{1-x^2}}dx
=\frac{1}{\pi}\int_0^\pi f(\cos \theta)d\theta
\]
and the characters, in the present identification, 
are given by $\psi_n(\cos\theta)=\cos n\theta$ for $n\iin\En_0$.

\begin{teespfa}
$\mathrm{Z}_{\Zee_2}\blone{\Tee}$ is amenable.
\end{teespfa}

\proof  We let $K_n:[0,\pi]\to\Ree^{\geq 0}$ for $n\in\En_0$ denote the well-known 
Fejer kernel (see \cite[2.5]{katznelson}, for example), so
\[
K_n=\sum_{k=0}^n\left(1-\frac{2k}{n+1}\right)\psi_k\comp\cos.
\]
Then let $\mu_n=\sum_{k=0}^n\left(1-\frac{2k}{n+1}\right)\psi_k\otimes\psi_k$.
We have that $(\mu_n)$ is an approximate diagonal by
Theorem \ref{theo:diagonal}.  Moreover, $(\mu_n)$ is bounded since
\begin{align*}
\norm{\mu_n}_1&=\frac{1}{\pi^2}\int_0^\pi\int_0^\pi\left|
\sum_{k=0}^n\left(1-\frac{2k}{n+1}\right)\cos k\theta\cos k\theta'\right|d\theta d\theta' \\
&=\frac{1}{2\pi^2}\int_0^\pi\int_0^\pi\left|
\sum_{k=0}^n\left(1-\frac{2k}{n+1}\right)\bigl(\cos k(\theta+\theta')+
\cos k(\theta-\theta')\bigr)\right|d\theta d\theta' \\
&=\frac{1}{2\pi^2}\int_0^\pi\int_0^\pi \bigl(K_n(\theta+\theta')+K_n(\theta-\theta')\bigr)
d\theta d\theta' =1.
\end{align*}
Thus $\mathrm{Z}_{\Zee_2}\blone{\Tee}\cong\blonek$ is amenable.  \endpf


\section{Some non-compact groups}

\subsection{Preliminaries and Notation}\label{ssec:pandn2}
If $G$ is a locally compact group, then $\zbloneg\not=\{0\}$ if and only
if $G$ has a relatively compact neigbourhood which is invariant under inner
automorphisms, i.e.\ $G$ is an $[IN]$-group; see \cite[Prop.\ 1]{mosak}.  
In fact, it is shown in \cite[Cor.\ 1.5]{liukkonenm}
that $\zbloneg$ is related to certain centers of $[FIA]^-_B$-groups,
which we define below.

Let $\autg$ denote the space of continuous automorphism of $G$
which can be endowed with a Hausdorff topology \cite[(26.5)]{hewittrI}.
We let $\inng=\{s\mapsto tst^{-1}:t\in G\}$ denote the group of inner
automorphisms in $\autg$.  We say $G$ has {\it relatively compact
inner automorphisms} if $\inng$ is realtively compact in $\autg$.
More generally, if there is a relatively compact subgroup $B$ of
$\autg$ such that $B\supset\inng$ we say $G$ is of class $[FIA]^-_B$.
We let for $\beta\iin B$ and $f\iin\bloneg$, $f\comp\beta(s)=f(\beta(s))$ for
almost every $s\iin G$.  We then let
\[
\zbbloneg=\{f\in\bloneg:f\comp\beta=f\text{ for all }\beta\iin B\}.
\]
This is a subalgebra of $\zbloneg$.  The 
result \cite[Cor.\ 1.5]{liukkonenm}, to which we alluded, above,
is that for an $[IN]$-group $G$, there is 
open normal subgroup $G_0$ of $G$ generated by all elements
with relatively compact conjugacy classes,
and a closed normal subgroup of $G_0$, $N$, which is the intersection
if all $\inng|_{G_0}$-invariant neighbourhoods of $e$, 
so that group $B=\{sN\mapsto t^{-1}stN:t\in G_0\}$ is relatively comapct
in $\auth$ where $H=G_0/N$, and
\begin{equation}\label{eq:liukkonenmosak}
\zbloneg\cong\zbbloneh.
\end{equation}
We let $\fX_B(G)$ denote the Gelfand spectrum of $\zbbloneg$, 
and let $\fX(G)=\fX_{\inng}(G)$.
The identification (\ref{eq:liukkonenmosak}) gives a natural identification
$\fX(G)\cong\fX_B(H)$.  It follows from \cite[4.12]{mosak0} 
(see \cite[4.2]{hulanicki}) that $\fX_B(G)$ may be 
identified with a certain family of continuous positive definite functions 
on $G$.

We record the following important structural result, which will be key to many
of the results which follow.  It summarises results from
\cite[Prop.\ 2.3]{liukkonenm} and \cite[Lem.\ 1]{mosak}.
See the summary presented in \cite{stegmeir}.

\begin{stegmeir}\label{lem:stegmeir}
Let $G$ be an $[FIA]^-_B$-group and suppose there exists
a compact $B$-invariant
subgroup $K$ such that each ``$\beta$-commutator'' $s^{-1}\beta(s)\in K$,
where $\beta\iin B$ and $s\iin G$ (thus $G/K$ is abelian).
Define an equivence relation on $\fX_B(G)$ by 
\[
\chi\sim\ome\;\iff\;\chi|_K=\ome|_K.
\]
Let $[\chi]$ denote the equivalence class of $\chi$.  Then

{\bf (i)} there is a family of ideals $\{J(\chi):[\chi]\in\fX_B(G)/\sim\}$ such that
\[
J(\chi)\cap J(\ome)=\{0\}\iif\chi\not\sim\ome, \quad
\zbbloneg=\wbar{\bigoplus_{[\chi]\in\fX_B(G)/\sim}J(\chi)}
\]
\phantom{mmm}
\parbox[t]{4.5in}{and each $J(\chi)$ is isomorphic to $\blone{G(\chi)}$, where
$G(\chi)$ is an abelian group, isomorphic to a quotient of an open subgroup
of $G$ by $K$; and}

{\bf (ii)} $\{\chi|_K:[\chi]\in\fX_B(G)/\sim\}$ is an orthogonal family
in $\mathrm{L}^2(K)$.
\end{stegmeir}






Note that for such a compact subgroup as $K$ to exist, it is necessary
and sufficient that the closed subgroup generated by $B$-commutators
be compact.  In this case $G$ is said to be an {\it $[FD]^-_B$-group}.
Note that if $G$ is compact we may take $K=G$ and we obtain,
for each $\pi\in\what{G}\cong\fX(G)$, $J(\chi_\pi)=\Cee\chi_\pi
\cong\blone{G/G}$.

\subsection{Some amenable centres}\label{ssec:fincommutator}
If $A,B$ are any pair of subgroups
of $G$, we let $[A,B]$ denote the closed subgroup generated
by commutators $\{aba^{-1}b^{-1}:a\in A,b\in B\}$.  The
derived subgroup is given by $G'=[G,G]$.

The following result is a generalisation of \cite[Thm.\ 1]{stegmeir}.
We recall that if $G$ is an $[IN]$-group, then the subgroup
$G_0$, of all elements with relatively compact conjugacy classes is an open normal
subgroup.

\begin{finitederived}\label{theo:finitederived}
If $[G,G_0]$ is finite, then $\zbloneg$ is amenable.
\end{finitederived}

\proof  We may suppose that $\zbloneg\not=\{0\}$, so $G$ has an invariant
neighbourhood.  Let $B$ and $H=G_0/N$ be as in (\ref{eq:liukkonenmosak}) and
$K=[G,G_0]/N$.  Then it is straighforward to check that $K$ is $B$-invariant
and that it is generated by $B$-commutators.
Since $K$ is finite, the orthogonality relations given in Lemma \ref{lem:stegmeir}  (ii)
imply that there are only 
finitely many ideals $\{J(\chi):[\chi]\in\fX_B(G)/\sim\}$.  It then follows
from Lemma \ref{lem:stegmeir}  (i), \cite[Prop.\ 5.2]{johnsonM}, 
and the fact that each $\blone{G(\chi)}$ is amenable,
that $\zbloneg\cong\zbbloneh$ is amenable.  \endpf

Observe that condition of the theorem above holds when $G'$ is finite.
It also holds when $G_0=\{e\}$, in which case $G$ is called an infinite conjugacy 
class group.





\subsection{Some hyper-Tauberian centres}\label{ssec:shytaualg}
We direct the reader to Section \ref{ssec:hypertauberian} for the
definition and consequences of the hyper-Tauberian property.

\begin{hypertauberian1}\label{prop:hypertauberian1}
Suppose $G$, $B$ and $K$ are as in the hypotheses of Lemma \ref{lem:stegmeir}.  
Then $\zbbloneg$ is hyper-Tauberian. 
\end{hypertauberian1}




\proof By Lemma \ref{lem:stegmeir} (i) we may write
\[
\zbbloneg=\wbar{\bigoplus_{[\chi]\in\fX_B(G)/\sim}J(\chi)}\cong
\wbar{\bigoplus_{[\chi]\in\fX_B(G)/\sim}\blone{G(\chi)}}.
\]
Hence it follows from \cite[Cor.\ 13]{samei}
that $\zbbloneg$ is hyper-Tauberian.  
 \endpf

We say that $G$ is an $[FC]^-$-group if each conjugacy class
in $G$ is relatively compact.  In the notation of Section \ref{ssec:pandn2}
this is the same as having $G=G_0$.

\begin{hypertauberian2}\label{theo:hypertauberian2}
If $G$ is an $[FC]^-$-group, then $\zbloneg$ is hyper-Tauberian.
\end{hypertauberian2}




\proof In the notation of (\ref{eq:liukkonenmosak}) we have that
$H=G/N$ and $B=\innh$.  Thus $\zbloneg\cong\zblone{H}$, and we
may assume $G$, itself, is an $[FIA]^-$-group.  

If $G$ is compactly generated, then \cite[(3.20)]{grosserm}
gaurantees that the derived group $K=G'$ is compact.  Hence we can
apply Proposition \ref{prop:hypertauberian1}, and we are done.

If $G$ is not compactly generated, we must localise our argument
to a compactly generated subgroup.  We first wish to see that
$\zcontcg$, the space of all of compactly supported continuous elements
of $\zbloneg$, is dense in $\zbloneg$.  We note that 
$P:\bloneg\to\zbloneg$, given for almost every $s\iin G$ by
$Pf(s)=\int_{\wbar{\inng}}f(\beta(s))d\beta$, defines a surjective quotient map.
Hence if $f\in\zbloneg$ and $(u_n)\subset\fC_c(G)$ is a sequence 
with $\lim_nu_n=f$, then $\lim_nPu_n=Pf=f$.

Now let $T:\zbloneg\to\zbloneg^*$ be a local operator.  To see that
$T$ is a $\zbloneg$-module map, it suffices to show that
\begin{equation}\label{eq:modulemap}
\dpair{T(u\con v)}{w}=\dpair{u\con T(v)}{w}
\end{equation}
for any $u,v,w\iin\zcontcg$.  The set $U=\{s\in G:|u(s)|+|v(s)|+|w(s)|>0\}$
is $\inng$-invariant, open and relativley compact.  Hence $U$ generates
a normal open subgroup $F$ of $G$.  We let $B=\inng|_F$ and note that
$F$ is an $[FIA]^-_B$-group.  

We have that the closed subgroup $K$ generated by $B$-commutators in $F$ 
is compact.  This is noted in \cite{lasser}, though does follow obviously from
\cite[(3.20)]{grosserm}.   Let us show how this can be proved
from \cite{grosserm}.  It is shown in \cite[(3.16)]{grosserm} that $G'$ consists 
of periodic elements, elements which individually generate relatively compact 
subgroups of $G$.  Hence $K=[F, G]\subset G'$ consists of periodic elements.
Since $F$ is compactly generated and an $[FIA]^-_B$-group, it is clear that
$K$ is compactly generated.  Then by \cite[(3.17)]{grosserm}, $K$ is compact.

Clearly $K$ is $B$-invariant.  Thus $\zbblonef$ is hyper-Tauberian by
Proposition \ref{prop:hypertauberian1}.  We note that $\zbblonef$ is the 
closed subalgebra of all elements of $\zbloneg$ which vanish almost everywhere
off of $F$.   Moreover, the mapping $\chi\mapsto\chi|_{\zbblonef}$ maps
$\fX(G)$ continuously onto $\fX_B(F)$, by \cite[Prop.\ 2.9]{liukkonenm}.
Let $\iota:\zbblonef\to\zbloneg$ be the injection map, so
$\iota^*\comp T\comp\iota:\zbblonef\to\zbblonef^*$ is a local map.
Then $\iota^*\comp T\comp\iota$ is a $\zbblonef$-module map.
Since $u,v,w\in\zbblonef$, we see that (\ref{eq:modulemap}) holds.  \endpf

We note that there are non-$[FC]^-$-groups for which the above result fails.
Let $n\geq 3$ and $G_n=\Ree^n\ltimes\mathrm{SO}(n)_d$, the semi-direct product
of $\Ree^n$ with the discrete special orthogonal group. 
We have for odd $n$ that $\zblone{G_n}\cong
\mathrm{Z}_{\mathrm{SO}(n)}\mathrm{L}^1(\Ree^n)$;
for $n=3$ this was observed in \cite[p.\ 162]{liukkonenm}.
(Note that for even $n$ we have $\z{\mathrm{SO}(n)}=\{1,-1\}=\Zee_2$
and we have 
$\zblone{G_n}\cong\mathrm{Z}_{\mathrm{SO}(n)}\mathrm{L}^1(\Ree^n\ltimes\Zee_2)$.)
It is proved in \cite[Prop.\ 2.6.8]{reiters} 
(see also \cite[Thm.\ 5.5]{azimifard})
that for $n\geq 3$, 
$\mathrm{Z}_{\mathrm{SO}(n)}\mathrm{L}^1(\Ree^n)$ admits non-zero
point derivations.  Hence this algebra cannot even be weakly amenable, 
neverless hyper-Tauberian, as noted in Theorem \ref{theo:hypertaubprop}.
Moreover, for $n\geq 3$,
it is shown \cite[2.6.10]{reiters} 
that except for the augmentation
character, no singleton in $\fX_{ \mathrm{SO}(n)}(\Ree^n)$ is a set of 
spectral 
synthesis.  




{
\bibliography{zlonegbib}
\bibliographystyle{plain}
}

\medskip
\noindent Ahmadreza Azimifard
\newline {\sc Fields Institute, 222 College Street, Toronto, Ontario, M5T\;3J1,
Canada}
\newline E-mail: {\tt aazminif@fields.utoronto.ca}

\medskip
\noindent Ebrhaim Samei:
\newline {\sc Department of Pure Mathematics, University of Waterloo,
Waterloo, Ontario, N2L\;3G1, Canada} 
\newline E-mail: {\tt esamei@uwaterloo.ca}

\medskip
\noindent Nico Spronk:
\newline {\sc Department of Pure Mathematics, University of Waterloo,
Waterloo, Ontario, N2L\;3G1, Canada} 
\newline E-mail: {\tt nspronk@uwaterloo.ca}

\end{document}